\documentclass[12pt]{article}
\usepackage{indentfirst,latexsym,bm,amsmath,amssymb,amsfonts,lscape,rawfonts}
\usepackage{algorithm,algorithmic}
\usepackage{graphicx}
\usepackage{cite,epic,eepic,euscript,verbatim,amsmath,amssymb,amsthm,afterpage,float,bm,paralist, amscd}
\usepackage{mathtools}
\usepackage{epsfig}
\usepackage{caption}
\usepackage{subcaption}
\usepackage{geometry}
\usepackage{MnSymbol}
\newtheorem{theorem}{Theorem}
\newtheorem{lemma}{Lemma}[section]
\newtheorem{proposition}{Proposition}[section]
\newtheorem{corollary}{Corollary}[section]
\newtheorem{remark}{Remark}
\newtheorem{definition}{Definition}

\numberwithin{figure}{section}

\newcommand{\bth}{\begin{theorem}}
\newcommand{\eth}{\end{theorem}}
\newcommand{\bpr}{\begin{proposition}}
\newcommand{\epr}{\end{proposition}}
\newcommand{\bde}{\begin{definition}}
\newcommand{\ede}{\end{definition}}
\newcommand{\blem}{\begin{lemma}}
\newcommand{\elem}{\end{lemma}}
\newcommand{\bco}{\begin{corollary}}
\newcommand{\eco}{\end{corollary}}
\newcommand{\prove}{\begin{proof}}
\newcommand{\done}{\end{proof}}
\newcommand{\brem}{\begin{remark}}
\newcommand{\erem}{\end{remark}}

\numberwithin{equation}{section}

\newcommand{\beq}{\begin{equation}}
\newcommand{\eeq}{\end{equation}}
\newcommand{\lb}{\label}
\newcommand{\eps}{\varepsilon}

\newcommand{\bbN}{{\mathbb{N}}}
\newcommand{\bbR}{{\mathbb{R}}}

\newcommand{\bbZ}{{\mathbb{Z}}}

\newcommand{\bbT}{{\mathbb{T}}}

\newcommand{\calH}{{\mathcal{H}}}
\newcommand{\calV}{{\mathcal{V}}}

\begin{document}

\title{Numerical Evidence of Exponential Mixing by Alternating Shear Flows}

\author{Li-Tien Cheng
\and
Frederick Rajasekaran
\and
Kin Yau James Wong
\and
Andrej Zlato\v s
\date{\today}
}

\newcommand{\Addresses}{{
  \bigskip
  \footnotesize

  \textsc{Department of Mathematics, UC San Diego, La Jolla, CA 92093, USA}\par\nopagebreak
  \textit{E-mail:} \texttt{l3cheng@ucsd.edu, zlatos@ucsd.edu}

  \medskip

  \textsc{UC San Diego, La Jolla, CA 92093, USA}\par\nopagebreak
  \textit{E-mail:} \texttt{frajasek@ucsd.edu, k1wong@ucsd.edu}
}}

%

\maketitle

\begin{abstract}
We performed a numerical study of the efficiency of mixing by alternating horizontal and vertical shear ``wedge'' flows on the two-dimensional torus.  Our results suggest that except in cases where each individual flow is applied for only a short time, these flows produce exponentially fast mixing.  The observed mixing rates are higher when the individual flow times are shorter (but not too short), and randomizing either the flow times or phase shifts of the flows does not appear to enhance mixing (again when the flow times are not too short).  In fact, the latter surprisingly seems to inhibit it slightly.
\end{abstract}

\section{Introduction and Motivation}

The study of mixing of substances by incompressible flows has practical applications
as well as connections to multiple branches of  mathematics and science.  The simplest mathematical model of the process of mixing in the absence of diffusion is the transport PDE
\beq \lb{1.1}
\rho_t+u\cdot\nabla\rho=0,
\eeq
where $\rho$ represents concentration of the mixed substance with some initial value $\rho(\cdot,0)=\rho_0$, and $u$ is the (prescribed and divergence-free) velocity of the mixing flow.  We will consider here $\rho\in L^\infty(\bbT^2\times[0,\infty))$, so the physical domain will be the two-dimensional torus.  Moreover,  since addition of a constant does not affect the dynamic of \eqref{1.1} and the spatial average $\int_{\bbT^2} \rho(x,t)dx= \int_{\bbT^2} \rho_0(x)dx$ of solutions is conserved, we can restrict our analysis to mean-zero solutions without loss.  We also stress that our interest is in mixing  by {\it divergence-free fluid flows} (due to practical considerations preferably time-periodic ones, possibly up to some simple transformations), acting on $\bbT^2$ continuously in time, as opposed to general {\it measure-preserving maps} $T:\bbT^2\to\bbT^2$ (which represent the discrete-time version of the problem but may not result from real-world advective stirring).

Two important questions about \eqref{1.1} concern optimal mixing rates of solutions, given certain natural constraints on the mixing flows $u$ (see the Section \ref{S2} for related definitions and further details), and which flows achieve these rates.   Addressing the first question, Crippa and De Lellis essentially showed in \cite{CL} that one cannot achieve faster than exponential-in-time mixing, thus also proving a modified version of Bressan's rearrangement cost  conjecture \cite{Bressan, B2}.  That this exponential rate is indeed achievable was shown by Yao and Zlato\v s \cite{YaoZla}, who found flows exponentially mixing any given initial data $\rho_0$ (including on domains with boundaries), as well as by Alberti, Crippa, and Mazzucato \cite{ACM2}, whose results only apply to a special class of initial data on $\bbT^2$ but also to a larger set of flow constraints.  These results therefore established optimality of exponential mixing and also found flows that achieve it.

Unfortunately, the flows constructed in \cite{YaoZla, ACM2} are quite complicated, far from time periodic, and heavily depend on the initial data.  All these facts have obvious practical limitations.  These issues were remedied by Elgindi and Zlato\v s \cite{ElgZla}, who constructed much simpler and time-periodic {\it almost universal exponential mixers} --- $\rho_0$-independent flows that mix exponentially all initial data that have at least some degree of regularity (the construction also extends to domains with boundaries and, unlike \cite{YaoZla, ACM2}, to all spatial dimensions).  These flows even mix all initial data asymptotically as $t\to\infty$, so they are {\it universal mixers}, but it was shown in \cite{ElgZla} that no universal mixer can have a  rate (exponential or otherwise) that is uniform in all bounded mean-zero $\rho_0$.

The construction in \cite{ElgZla} nevertheless still has one limitation.  While the constructed flows are H\" older continuous in space, they are not Lipschitz and their flow maps are discontinuous.  (The flows in \cite{YaoZla} share this limitation;  those in \cite{ACM2} apply to solutions taking only two values, so they can be modified arbitrarily on each of the two level sets without changing the solution dynamics, which allows one to avoid potential singularities in their construction.)  Also, Bedrossian, Blumenthal, and Punshon-Smith \cite{BBP} showed that generic solutions $u$ to the 2D Navier-Stokes equations with certain stochastic forcings are almost universal exponential mixers, which is of obvious practical interest.  On the other hand, these flows  are again quite complicated and not time-periodic, as well as not deterministic, and they only satisfy the required constraints on average in time rather than pointwise.  It is therefore still an open question whether (time-periodic) smooth or at least Lipschitz continuous almost universal exponential mixers (or even just universal mixers) exist on $\bbT^2$ or other domains.

One candidate for such flows on $\bbT^2$ was proposed by Pierrehumbert \cite{Pie,Pie2},
and this suggestion is quite simple although not time-periodic.  It is almost every realization of the random vector field taking values $(\sin(2\pi x_2+\omega_n),0)$ and $(0,b\sin(2\pi x_1+\omega_n))$ (with $b\in\bbR$ some constant) on time intervals $(n-1,n-\frac 12]$ and $(n-\frac 12, n]$ (for $n\in\bbN)$, respectively.  Here  $\omega_n$ are independent random variables uniformly distributed over $\bbT$; note also that while these flows are not continuous in time, this is easily remedied by a simple reparametrization described in \cite{YaoZla}.

These flows are a representative of a wider class of alternatively horizontal and vertical shear flows. Heuristically, they appear to have very good mixing properties in many cases, but we are not aware of any rigorous proofs.
The goal of the present work is to demonstrate numerically that such flows can indeed yield exponential mixing of passive scalars advected by them.  We will consider random flows, with randomness in phase and/or flow time (the latter will replace the amplitude $b$ above), as well as deterministic time-periodic ones.

One difficulty with a computational approach to \eqref{1.1} is that when mixing is fast, solutions quickly become very rough, which poses a challenge from the numerical standpoint.  This may be further amplified when one considers non-smooth initial data, which we will do here because
\beq \lb{1.2}
\rho_0:= \chi_{[0,1/2)\times[0,1)} - \chi_{[1/2,1)\times[0,1)}
\eeq  
(considered also in, e.g., \cite{Bressan, ACM2}) is a natural choice of a ``minimally premixed''  initial datum; of course, $\rho_0$ is smooth away from the line $x_1=\frac 12$.

This problem can be somewhat remedied by adding a smoothing diffusion term to \eqref{1.1}, 
but we do not wish to take this route and will instead address the issue by using a setup in which we can minimize the resulting complications.  We will use the approach from the advection step of Pierrehumbert's lattice method \cite{Pie2} (but with no diffusion step), where the solution is approximated by a linear combination of characteristic functions of $2^{2N}$ squares of size $2^{-N}\times 2^{-N}$ from a fine grid into which $\bbT^2$ is split (with $N\in\bbN$).  These squares are then moved according to the prescribed shear flow, with each shift rounded to be an integer multiple of the grid scale $2^{-N}$.  This results in a specific permutation of these squares in each advective step.  Of course, one can equivalently represent each square by its ``lower left'' vertex, and these vertices are then the grid points from
\beq \lb{1.5}
G_N:= \left\{ 0, \frac 1{2^N},\dots,\frac {2^N-1}{2^N} \right\}^2 \subseteq\bbT^2,
\eeq
 whose coordinates are integer multiples of $2^{-N}$.   We will do so and thus have $\rho(\cdot, t)\in L^\infty(G_N)$.

\begin{figure}[ht]
    \centering
    \begin{subfigure}[t]{0.3\textwidth}
        \centering
        \includegraphics[width=\textwidth]{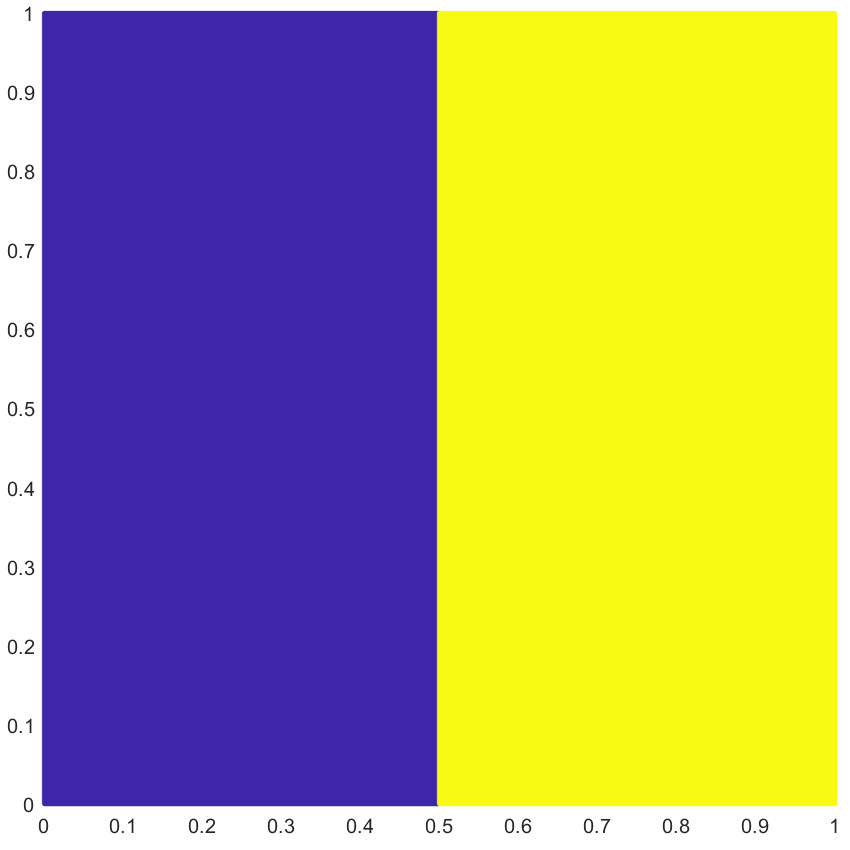}
        \caption{$\rho_0$}
    \end{subfigure}
        \hfill
    \begin{subfigure}[t]{0.3\textwidth}
        \centering
        \includegraphics[width=\textwidth]{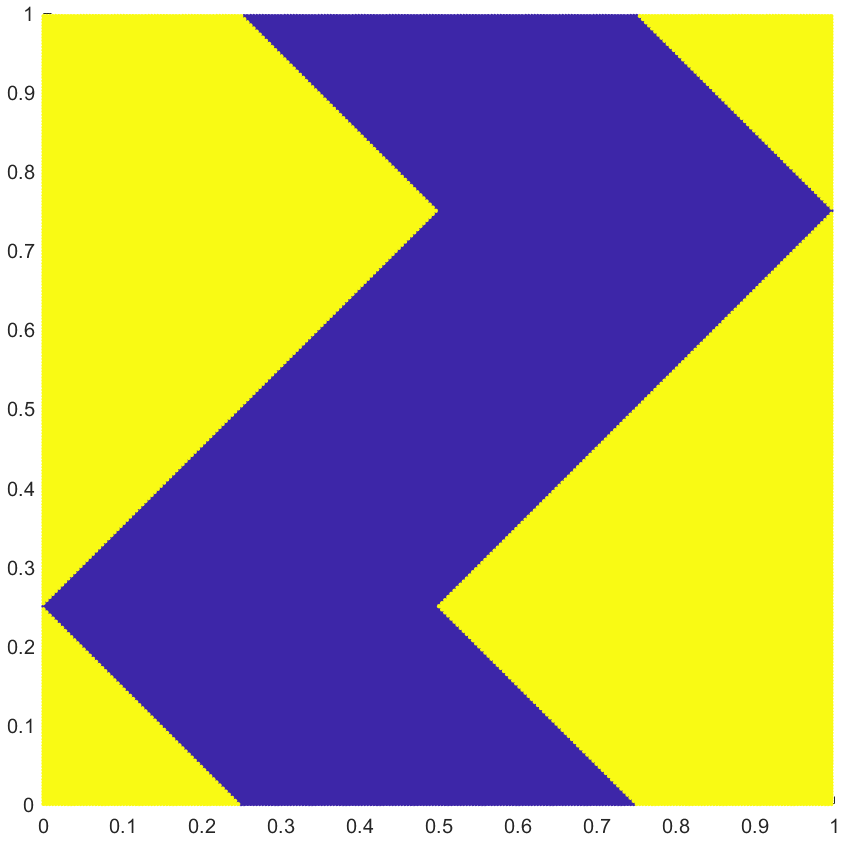}
        \caption{$\calH_{1/4}^{1}[\rho_0]$}
    \end{subfigure}
    \hfill
    \begin{subfigure}[t]{0.3\textwidth}
        \centering
        \includegraphics[width=\textwidth]{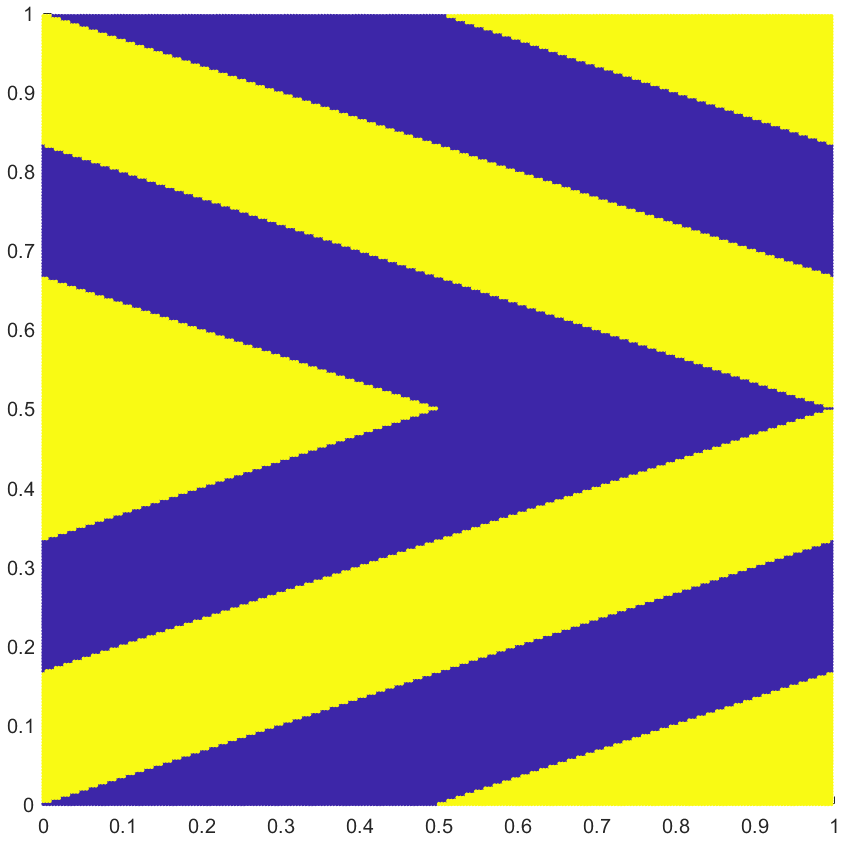}
        \caption{$\calH_{0}^{3}[\rho_0]$}
    \end{subfigure}

    \caption{Action of horizontal wedge flows on $\rho_0$ from \eqref{1.2} (values 1 and $-1$ are represented by blue and yellow colors, respectively).}
    \label{flow_examples}
\end{figure}

Moreover, we will avoid having to round the shifts by considering the
horizontal and vertical  ``wedge'' flows
\beq \lb{1.3}
    u^H_{\omega}\left(    x_1,x_2\right) := \left( d_{\mathbb{T}}(x_2, \omega), 0 \right) \qquad\text{and}\qquad
    u^V_{\omega}\left(    x_1,x_2\right) := \left( 0, d_{\mathbb{T}}(x_1, \omega) \right)
 \eeq
instead of the sine flows from \cite{Pie}, with 
\[
d_{\mathbb{T}}(x,y) := \min \left\{ |x - y|, 1 - |x - y| \right\}
\]
 the distance on $\bbT$ (so $d_{\mathbb{T}}(x,y)\in [0,\frac 12 ]$) and $\omega\in\bbT$ some {\it phase shift}.  If $\omega$ is an integer multiple of $2^{-N}$ and either of these flows is applied for an integer length of time $\tau\in\bbZ$, then the resulting time-$\tau$ flow maps
\beq \lb{1.4}
    H_{\omega}^{\tau}\left(    x_1,x_2\right) := \left(x_1+ \tau d_{\mathbb{T}}(x_2, \omega), x_2\right) \qquad\text{and}\qquad
    V_{\omega}^{\tau}\left(x_1,x_2\right) :=\left(x_1, x_2+ \tau d_{\mathbb{T}}(x_1, \omega) \right)
\eeq
acting on $\bbT^2$ keep $G_N$  invariant. 
These then transform functions $f\in L^\infty(\bbT^2)$ via
\beq \lb{1.6}
    \calH_{\omega}^{\tau}[f] \left(    x_1,x_2\right) := f\left(x_1- \tau d_{\mathbb{T}}(x_2, \omega), x_2\right) \qquad\text{and}\qquad
    \calV_{\omega}^{\tau}[f] \left(x_1,x_2\right) := f\left(x_1, x_2- \tau d_{\mathbb{T}}(x_1, \omega) \right).
\eeq
Since we want our flows to satisfy time-uniform constraints,  $\tau\in\bbZ$ will represent the length of time  during which the particular shear flow is acting rather than the flow amplitude, and we will refer to it as \textit{flow time}.  For instance, the solution to \eqref{1.1} with $u=u_\omega^H$ is given by $\rho(\cdot,\tau)=\calH_{\omega}^{\tau}[\rho(\cdot,0)]$ (see Figure~\ref{flow_examples}).  We note that while these flows are only Lipschitz, which nevertheless still guarantees continuity of their flow maps, this comes with an additional advantage over the sine flows from \cite{Pie} and smooth flows in general.  Since they do not have ``flat'' spots (such as at the maxima and minima of the sines where their derivatives vanish), where less mixing is happening due to much less advective stretching in those regions, they should be much better candidates for efficient mixers.

\begin{figure}[htbp]
    \centering
    \begin{subfigure}[t]{0.3\textwidth}
        \centering
        \includegraphics[width=\textwidth]{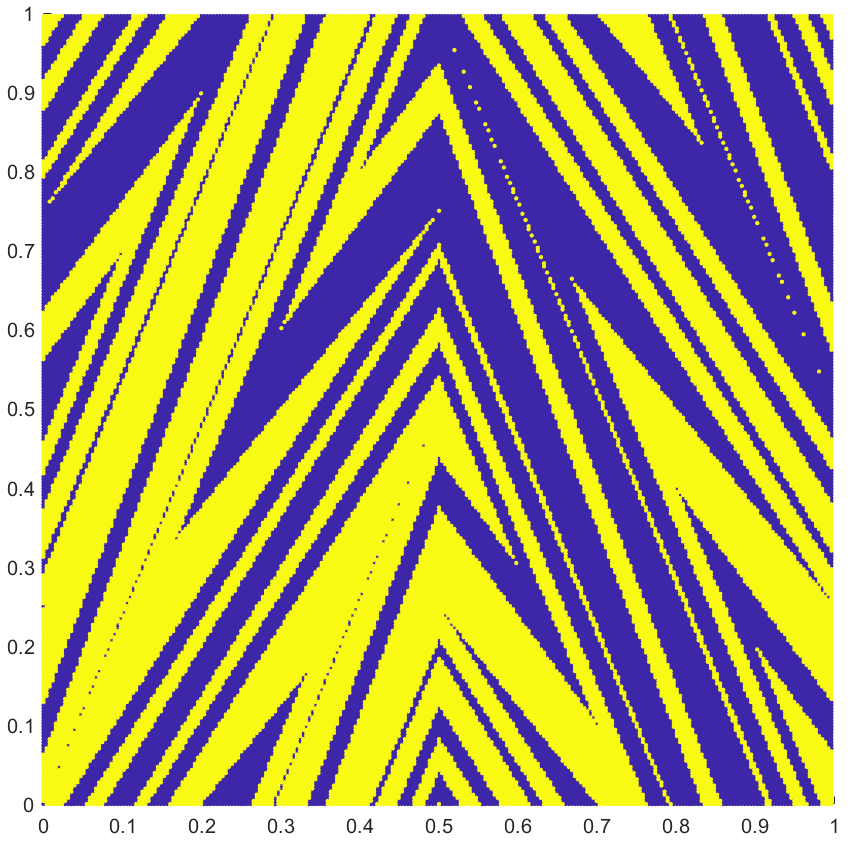}
        \caption{$(\calV_{0}^{2} \circ \calH_{0}^{2})^2[\rho_0]$}
    \end{subfigure}
    \hfill
    \begin{subfigure}[t]{0.3\textwidth}
        \centering
        \includegraphics[width=\textwidth]{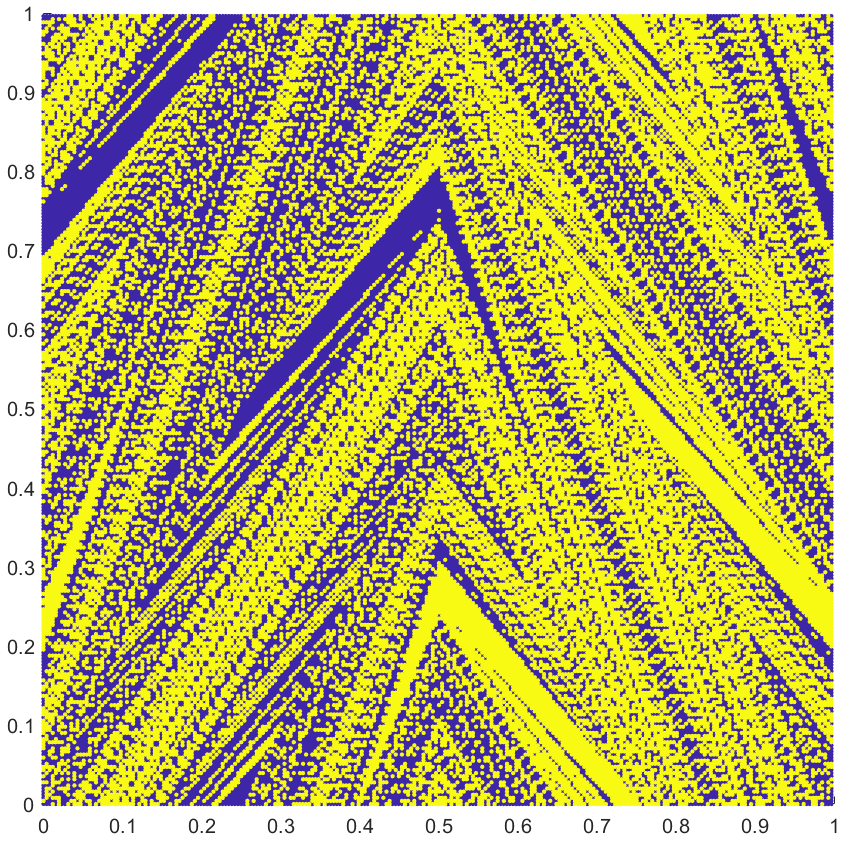}      
        \caption{$(\calV_{0}^{2} \circ \calH_{0}^{2})^{4}[\rho_0]$}
    \end{subfigure}
    \hfill
    \begin{subfigure}[t]{0.3\textwidth}
        \centering
        \includegraphics[width=\textwidth]{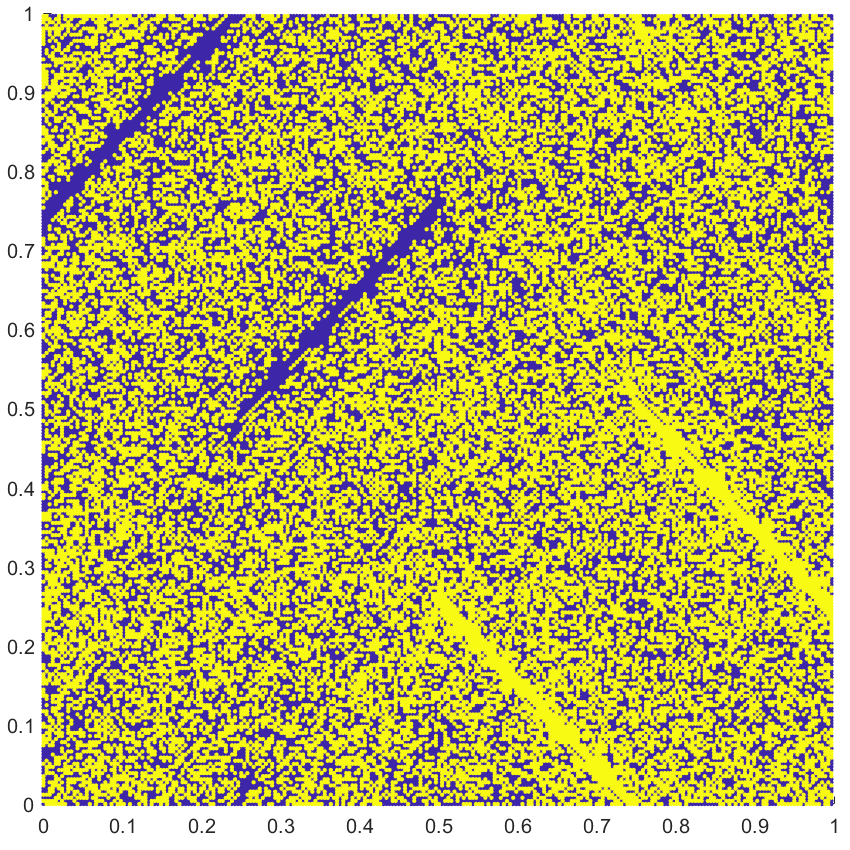}
        \caption{$(\calV_{0}^{2} \circ \calH_{0}^{2})^{8}[\rho_0]$}
    \end{subfigure}
    
    \caption{Action of alternating wedge flows with fixed $\tau$ and $\omega$ (the monochromatic  lines in the top left and bottom right of the pictures are due to $\tau=2$, see Subsection \ref{S3.2}).}
    \label{mixing_example}
\end{figure}

The above setup means that we will consider here functions $\rho_0,\rho_1,\dots\in L^\infty(G_N)$, with $\rho_{k+1}$  being either $\calH_\omega^1[\rho_k]$ or $\calV_\omega^1[\rho_k]$, where $\omega\in\{0,\frac 1{2^{N}},\cdots,\frac{2^N-1}{2^{N}}\}$ is the phase shift of the either horizontal or vertical wedge flow that acts during the integer-length time interval that contains $[k,k+1)$.
We will study four cases here, with both $\omega$ and $\tau$ fixed as well as random.  When both are fixed, the resulting  flow is $2\tau$-periodic in time (one could fix separate values of $\tau$ for the horizontal and vertical flows but we will not do this here).  In the random case, we will randomly choose new $\omega$ and/or $\tau$ each time we switch the direction of the flow between horizontal and vertical (see Section \ref{S2} below for details).  In Figure~\ref{mixing_example}  we demonstrate the action of alternating wedge flows, with both phase shift and flow time fixed,  on the initial datum $\rho_0$ from \eqref{1.2}. 
In  the rest of this paper, we restrict $\rho$ and $\rho_0$ to $G_N\subseteq \bbT^2$ with $N=15$, where each point from $G_N$ represents its adjacent $2^{-N}\times 2^{-N}$ grid square.

Our results support the conjecture that all the flows studied here are indeed exponential mixers for $\rho_0$,
except possibly some of those with $|\tau|\le 2$.  Although we only consider the single initial datum \eqref{1.2}, the fact that it is ``minimally premixed'' and that it appears to be  mixed exponentially quickly regardless of the choice of the phase shifts and flow times, suggest that  results for other initial data would be similar, and therefore these flows could in fact be almost universal exponential mixers.

Let us now turn to the specifics of our work.  We discuss the details of our setup and the definition of mixing scales that we will use here in Section \ref{S2}.  We then present our results in Section \ref{S3}, the main one being numerical evidence of exponential mixing by alternating wedge flows, summarized in the table in Figure \ref{main_data}.
\medskip

{\bf Acknowledgements}.  AZ thanks Gautam Iyer and Jean-Luc Thiffeault for illuminating discussions.   LTC acknowledges partial support by  NSF grant DMS-1913144.  FR and KYJW were supported in part   by Division of Physical Sciences Undergraduate Summer Research Awards and 
TRELS Awards from UC San Diego.
AZ acknowledges partial support by  NSF grant DMS-1900943 and by a Simons Fellowship.

\section{Definition of Mixing Scales and the Modeled Flows} \lb{S2}

\subsection{Mixing Scales of Bounded Functions}

In order to be able to study mixing efficiency of flows, we need to adapt a relevant definition of mixing scales of solutions $\rho$ to \eqref{1.1} with initial data $\rho_0$.  Following \cite{Bressan,YaoZla}, one option is to say that $\rho(\cdot, t)$ is $\kappa$-mixed to scale $\eps>0$  (for some $\kappa\in(0,1)$) when
\beq \lb{2.1}
\left| \strokedint_{B_{\varepsilon}(y)} \rho(x,t)\,dx \right| \leq \kappa ||\rho_0||_{\infty}
\eeq
holds for each $y\in\bbT^2$ (note that $||\rho||_{\infty}=||\rho_0||_{\infty}$).  The  mixing scale of $\rho(\cdot,t)$ is then the infimum of all $\eps> 0$ such that $\rho(\cdot,t)$ is mixed to scale $\eps$.  This is also called the {\it geometric mixing scale} in \cite{YaoZla}, and was recently used in various other works including \cite{ACM, ACM2,CL, ElgZla, CLS}.  Other alternatives are the {\it functional mixing scale} $\|\rho(\cdot,t)\|_{\dot H^{-s}}\|\rho(\cdot,t)\|_{\infty}^{-1}$ for some $s>0$ (particularly $s\in\{\frac 12,1\}$), used for instance in  \cite{ACM2, BBP, IKX, LTD, LLNMD, S,CLS, ElgZla} and closely related to the geometric mixing scale (see \cite{YaoZla}), as well as the Wasserstein distance of the positive and negative parts of $\rho(\cdot,t)$ \cite{BOS, OSS, S, S2}.
We note that there is also a large literature on the interaction between mixing and diffusion, as well as alternative definitions of mixing in the diffusive setting, but we will not attempt to provide an overview here and instead refer the reader to the review \cite{T} (which also concerns the diffusion-less case \eqref{1.1}) and references therein.

We will consider here a version of the geometric mixing scale above, which is most suited to our setting, but with further adjustments.  In the discrete setting of the grid $G_N$ that we will consider here, averaging the solution over discs may be problematic as the number of points from $G_N$ inside the disc $B_\eps(y)$ is not a constant multiple of $\eps^2$.  Moreover, finding averages over all discs of a particular radius centered at points from $G_N$ is unnecessarily computationally intensive, while restricting this to only points from some sub-grid would mean that not all points from $G_N$ are equally represented in the mixing scale computation.  

It is therefore both more reasonable and better suited to our setup to consider only $\eps=2^{-n}$ for $n=0,1,\dots$ and replace the discs $B_\eps(y)$ in \eqref{2.1} by all $2^{2n}$ squares
\[
S_n^{i,j} := \left[\frac{i}{2^n},\frac{i+1}{2^n}\right)\times\left[\frac{j}{2^n},\frac{j+1}{2^n}\right)
\]  
with $i,j\in \{0,1,\dots,2^{n}-1\}$.
One can easily show (see the proof of Lemma 2.7 in \cite{ElgZla}) that the resulting definition of the mixing scale (which will always be a power of $\frac 12$) is essentially equivalent to the above definition of geometric mixing scale when it comes to the study of asymptotic mixing rates (ratio of one with any $\kappa\in(0,1)$ and the other with any $\kappa'<\kappa$ is bounded above by a constant depending only on $\kappa-\kappa'$).  As mentioned in the introduction, we will also have $\rho_0,\rho_1,\dots\in L^\infty(G_N)$ instead of $\rho\in L^\infty(\bbT^2\times[0,\infty))$, so the averages over these squares will be just the averages over the points from $G_N$ contained in them.  Of course, this means that  the minimal possible mixing scale will be $2^{1-N}$ (unless $f\equiv 0$).
The choice of $\kappa\in(0,1)$ does not affect the exponential mixing rates as $t\to\infty$ (in the continuous space setting as well as in the $N\gg t$ regime in the discrete setting) but will have some effect on finite time intervals.  We choose $\kappa:=\frac 13$ as in \cite{Bressan}, which finally yields the following definition.

\bde \lb{D.2.1}
We say that a mean-zero function $f:G_N\to\bbR$ is {\it mixed to scale} $2^{-n}$ for some $n\in\{0,1,\dots,N\}$ if for each pair $i,j\in \{0,1,\dots,2^{n}-1\}$ we have
\[
\left| {2^{-2(N-n)}}\sum_{x\in\mathcal{S}_n^{i,j}\cap G_N} f(x)\right| \leq \frac{\| f\|_\infty}{3}.
\]
The {\it mixing scale} of $f$ is the smallest such $2^{-n}$.
\ede

Note that in our setting we will always have $\|\rho_k\|_\infty=\|\rho_0\|_\infty = 1$ for all $k\in \bbN_0$.

\subsection{The Modeled Flow Types}
\label{S2.2}

As discussed in the introduction, we will consider here four basic flow types, all alternating wedge flows.  Two will have fixed phase shifts (chosen randomly at the start) and two will have their phase shifts chosen randomly each time we switch the flow direction.  Two will have fixed flow times and two will have their  flow times chosen randomly each time we switch the flow direction. We will use the discretized framework described in the introduction with $N=15$ (so grid scale will be $2^{-15}$), modeling the flow dynamic via the mappings from \eqref{1.6} applied to the initial data \eqref{1.2} on the domain \eqref{1.5}.  The specific details are as follows.

\medskip\noindent
\textbf{Fixed Shift Fixed Time (FSFT):} 
We choose randomly phase shifts $\omega,\omega'\in\{0,\frac 1{2^{N}},\cdots,\frac{2^N-1}{2^{N}}\}$ (with uniform joint distribution), then fix these and some flow time $\tau\in\{2,\dots,10\}$, and  let 
\[
\rho_{k+1}:=
\begin{cases}
\calH_{\omega}[\rho_k] & k\in [2j\tau,(2j+1)\tau) \text{ for some $j\in\bbN_0$}, 
\\ \calV_{\omega'}[\rho_k] & k\in [(2j+1)\tau,(2j+2)\tau) \text{ for some $j\in\bbN_0$}, 
\end{cases}
\]
for $k=0,1,\dots$.
 Since one should expect different behavior for different  $\tau$ (which we do confirm), we model these cases  separately.  The phase shifts $\omega,\omega'$ are not expected to have a significant effect on the mixing rates (which we also confirm), although they will have some effect on the computed mixing scales at individual times $k$.  This is clear for the horizontal shift $\omega'$ of the vertical  flow $u^V_{\omega'}$ because different shifts align differently with $\rho_0$.  The vertical shift $\omega$ would have no effect on the mixing scales if we were to  include in Definition~\ref{D.2.1} all the $2^{-n}\times 2^{-n}$ squares with vertices in $G_N$ rather than just the squares $S_n^{i,j}$ (which have vertices in $G_n$).  We do not do this in order to shorten the required computing time, and our simulations show that the effect on the obtained mixing rates would also be negligible.  We run the simulation 100 times (i.e., with 100 random choices of $(\omega,\omega')$) for each $\tau$.
 
\medskip\noindent
\textbf{Random Shift Fixed Time (RSFT):} 
Here the phase shifts are i.i.d.~random variables $\omega_0,\omega_0', \omega_1,\omega_1',\dots \in\{0,\frac 1{2^{N}},\cdots,\frac{2^N-1}{2^{N}}\}$ (with uniform distribution) and the flow time $\tau\in\{2,\dots,10\}$ is again fixed, so we have
\[
\rho_{k+1}:=
\begin{cases}
\calH_{\omega_j}[\rho_k] & k\in [2j\tau,(2j+1)\tau) \text{ for some $j\in\bbN_0$}, 
\\ \calV_{\omega_j'}[\rho_k] & k\in [(2j+1)\tau,(2j+2)\tau) \text{ for some $j\in\bbN_0$}, 
\end{cases}
\]
for $k=0,1,\dots$.  Contrasting the results in this case with those for FSFT allows one to see whether in the latter case the mappings $\calV_{\omega'}\circ\calH_{\omega}$ (which generate the FSFT dynamic, and coincide for all $(\omega,\omega')$ up to translation) may involve structures that slow down or accelerate mixing, since such structures would not persist in the RSFT case. We perform 100 simulations for each $\tau$.
 
\medskip\noindent
\textbf{Fixed Shift Random Time (FSRT):} 
Here the phase shifts are randomly chosen at the start and then kept fixed as in FSFT, but the flow time  is chosen randomly at each direction switch to see whether this can improve mixing. Since our FSFT and RSFT results show that the mixing rate depends nontrivially on the flow time when the latter is fixed, decreasing as $\tau$ increases from 3 or 4 to higher values, it makes sense to limit the randomness in the flow time to small intervals.  We therefore let the flow times be i.i.d.~random variables $\tau_0,\tau_0',\tau_1,\tau_1',\dots\in\{\tau-1,\tau,\tau+1\}$ (with uniform distribution) for some fixed $\tau\in\{2,\dots,10\}$, so
\[
\rho_{k+1}:=
\begin{cases}
\calH_{\omega}[\rho_k] & k\in [t_j,t_j+\tau_j) \text{ for some $j\in\bbN_0$}, 
\\ \calV_{\omega'}[\rho_k] & k\in [t_j+\tau_j,t_{j+1}) \text{ for some $j\in\bbN_0$}, 
\end{cases}
\]
for $k=0,1,\dots$, with $t_j:=\sum_{l=0}^{j-1}(\tau_l+\tau_l')$.  We  perform 100 simulations for each $\tau$.

\medskip\noindent
\textbf{Random Shift Random Time (RSRT):} 
Here the phase shifts $\omega_0,\omega_0', \omega_1,\omega_1',\dots$ are chosen as in RSFT and flow times $\tau_0,\tau_0',\tau_1,\tau_1',\dots$ as in FSRT, so now
\[
\rho_{k+1}:=
\begin{cases}
\calH_{\omega_j}[\rho_k] & k\in [t_j,t_j+\tau_j) \text{ for some $j\in\bbN_0$}, 
\\ \calV_{\omega_j'}[\rho_k] & k\in [t_j+\tau_j,t_{j+1}) \text{ for some $j\in\bbN_0$} ,
\end{cases}
\]
for $k=0,1,\dots$, with $t_j:=\sum_{l=0}^{j-1}(\tau_l+\tau_l')$.  We  perform 100 simulations for each $\tau\in\{2,\dots,10\}$.

\medskip
Of course, the formula for $\rho_{k+1}$ in the RSRT case also applies in the other cases, but with $(\omega_j,\omega_j'):=(\omega,\omega')$ for all $j\in\bbN_0$ and/or  $(\tau_j,\tau_j'):=(\tau,\tau)$ for all $j\in\bbN_0$.

\subsection{Computation of Exponential Mixing Rates} \lb{S2.3}

In each of the four  cases above and for each $\tau\in\{2,\dots,10\}$, we performed 100 simulations.  In every simulation we found the mixing scale $2^{-n_k}$ of the solution $\rho_k$ at each time $k\in\bbN_0$ via Definition \ref{D.2.1}.
We then used these mixing scales to find an exponential mixing rate of the flow in each individual run via linear regression over a relevant time interval (see below) and finally averaged these rates over the 100 runs.

\begin{figure}[htbp]
    \centering
        \includegraphics[scale = 0.55]{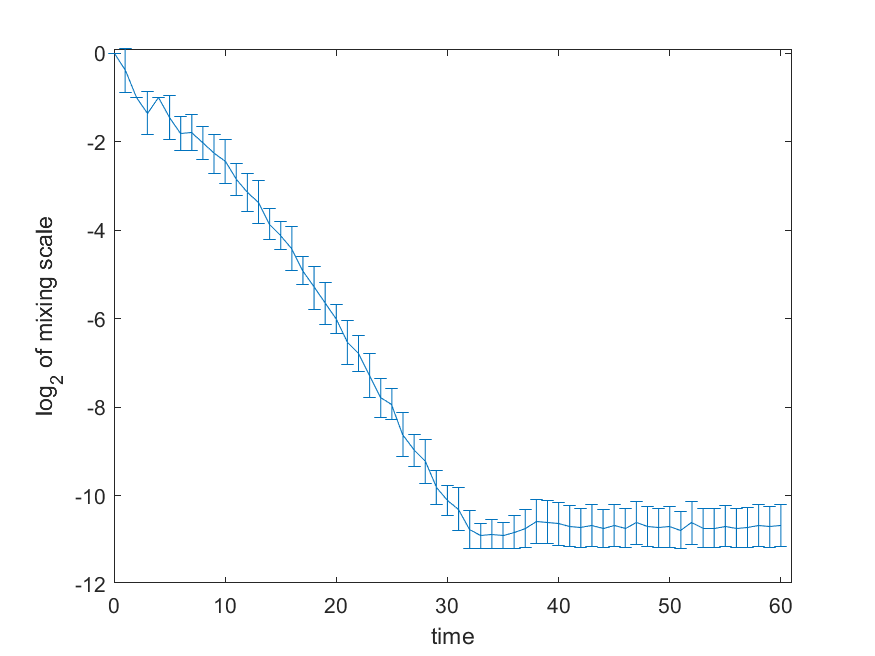} 
        \includegraphics[scale = 0.55]{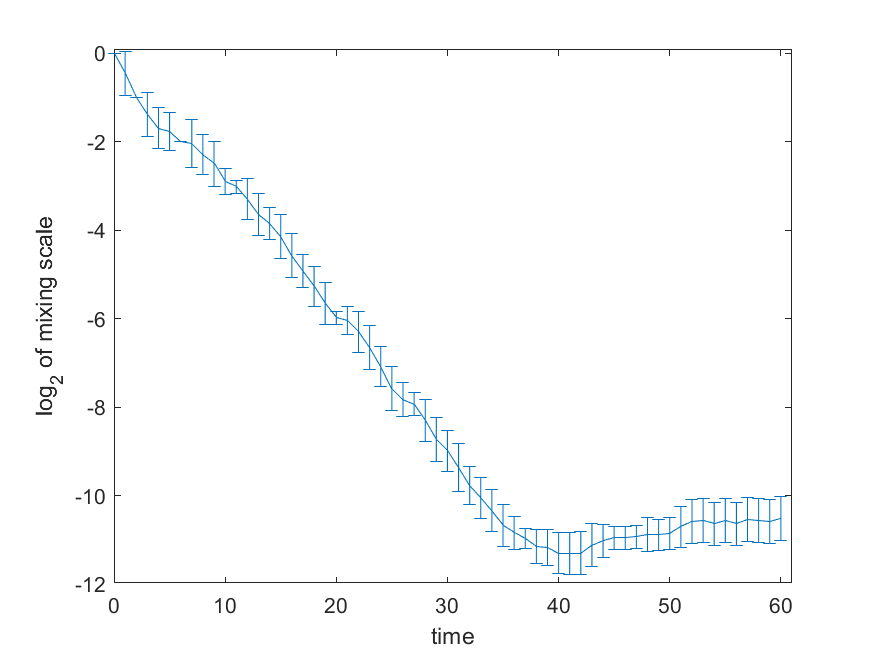}        
    \caption{Means and standard deviation error bars of the binary logarithms of the mixing scales for the 100 simulations 
   with $N = 15$ in the RSFT case with $\tau=3$ (left) and  in the FSRT case with $\tau=7$ (right).  The (mean) mixing scales have a near-uniform exponential decay before eventually plateauing around $2^{-11}$.}
    \label{plot_example}
\end{figure}

The computed mixing scales can never reach the grid scale $2^{-N}$, and they typically plateaued around $2^{4-N}$ in our simulations  with $N=15$ (see Figure \ref{plot_example}) as well as for other values of $N$ (since the square $S_{N-3}^{i,j}$ has 64 points from $G_N$, reaching mixing scale $2^{3-N}$ requires it to contain between 22 and 42 points of either color {\it for each} $(i,j)$).  In order to suppress this grid scale effect, we chose the end of the time interval for computation of the mixing rate to be the first time when the mixing scale reached $2^{5-N}$.  For $N=15$ this is $2^{-10}$, and we denote this time $T_{10}$ below.

Our simulations also showed that there is an initial time interval where the mixing scale may display somewhat irregular behavior.  An example of this is  in Figure \ref{plot_example}, which contains means and standard deviation error bars of the binary logarithms $-n_k$ of the mixing scales at different times $k$  for the 100 simulations  in two flow cases, the RSFT case with $\tau=3$ and the FSRT case with $\tau=7$.  One can observe near-plateaus of the averaged mixing scales in the time interval $[3,7]$ (roughly while these scales are between $2^{-1}$ and $2^{-2}$), likely due to an interaction of the initial data with different phase shifts, before they start an almost constant rate descent (until they plateau around $2^{-11}$).  A similar pattern appears before time $8$ in most of the other cases of flows with $\tau\ge 3$ (which are the ones providing efficient mixing, see Section \ref{S3} below), more so for smaller values of $\tau$.  None is more pronounced than on the left of Figure \ref{plot_example}, and therefore their effect on the computed mixing rates would be very small.  Nevertheless, in the interest of obtaining the most accurate approximations of the actual asymptotic mixing rates, we suppressed it by choosing the start of the time interval for computation of the mixing rates to be 8.  (We note that in almost all cases with $\tau\ge 3$, the first time when the averaged mixing scales dropped below $2^{-2}$ was either 7 or 8.  The exceptions were FSRT and RSRT with $\tau=3$, when this time was 9; these are however also somewhat exceptional, see Subection \ref{S3.1} below.)


So with $N=15$, we considered the time interval $[8,T_{10}]$ for each individual simulation (with $T_{10}$ being simulation-dependent).  We then found the corresponding exponential mixing rate 
as the {\it negative} of the slope of the line that best (in terms of least squares) fits
 the computed  binary logarithms $-n_k$ of mixing scales at all the integer times within this time interval (so this is then a base-2 exponential rate).  
We also computed the R-squared value of the fitting line (its proximity to 1 indicates a good fit and hence a near-uniform exponential decay of the mixing scales), as well as averages and standard deviations of the mixing rates and the R-values for the 100 runs for each flow type and each $\tau$.

\section{Main Results and Discussion} \label{S3}

In this section, we first present the results of the simulations described in Section \ref{S2}.  These constitute numerical evidence that the alternating wedge flows from Subection \ref{S2.2} indeed exhibit exponential mixing, except in some cases with small $\tau$, and we also identify the apparently most efficient mixers among these flows.  We then discuss the cases with small $\tau$, when mixing rates are lower and mixing can even be algebraic, as well as the related existence of structures fixed by the mixing dynamic.

\subsection{Evidence of Exponential Mixing} \label{S3.1}

The table in Figure \ref{main_data} contains averages of the  base-2 exponential mixing rates  for all four flow types and all $\tau\in\{2,\dots,10\}$, computed as described in Subsection \ref{S2.3} (with $N=15$ and 100 runs in each case).  It suggests that the most efficient mixing by the alternating wedge flows considered here happens when each individual wedge flow acts for the same time, which is either 3 or 4, and the phase shifts of the flows are not varied during each run.

It also shows that flows with $\tau=2$ are much worse mixers than those with $\tau\ge 3$; this is even more pronounced for $\tau=1$, which is why we did not include that case in our simulations.  We discuss these issues in Subsection \ref{S3.2} below, so will now concentrate on the cases with $\tau\ge 3$.
We note that  in all 32 of them, the means of the R-squared values were at least $0.9660$ and their
 standard deviations were no more than 0.0109, which is why we do not report these here (we note that even the 3200 individual R-squared values were all no less than 0.9290).  These numbers indicate an excellent fit and near-uniform exponential decay in all cases.

\begin{figure}[htbp]
\centering    
\begin{tabular}{|c|c|c|c|c|}
	\hline
	$\tau$ & FSFT & RSFT & FSRT & RSRT \\\hline\hline
	2 & --- & 0.2137 & 0.1279 & 0.1220 \\\hline
	3 & 0.3954 & 0.3781 & 0.3560 & 0.3480 \\\hline
	4 & 0.3955 & 0.3726 & 0.3723 & 0.3595 \\\hline
	5 & 0.3542 & 0.3342 & 0.3560 & 0.3429 \\\hline
	6 & 0.3409 & 0.3258 & 0.3293 & 0.3195 \\\hline
	7 & 0.3121 & 0.3046 & 0.3125 & 0.3070 \\\hline
	8 & 0.2959 & 0.2891 & 0.2946 & 0.2861 \\\hline
	9 & 0.2789 & 0.2723 & 0.2788 & 0.2714 \\\hline
	10 & 0.2631 & 0.2595 & 0.2620 & 0.2556 \\\hline
\end{tabular}
    \caption{Mean  (base-2) exponential mixing rates.}
    \label{main_data}
\end{figure}

For flow types FSFT and RSFT, the largest mixing rates were obtained for $\tau=3$ and $\tau=4$; those for $\tau\ge 5$ are noticeably lower, and steadily decrease as $\tau$ grows.  This phenomenon can be explained by noticing that mixing scales for time-independent shear flows decrease no faster than $O(t^{-1})$, and exponential mixing therefore results from alternation of horizontal and vertical flows.  This alternation cannot be too fast (as the cases $\tau=1,2$ show), but once each individual flow acts for a long enough time (which our simulations suggest to be 3 or 4), switching the flow direction becomes more beneficial for fast mixing than keeping it.  The reason for smaller mixing rates in the cases FSRT and RSRT with $\tau=3$ vs.~$\tau=4$ is the fact that the set  $\{\tau-1,\tau,\tau+1\}$, from which we randomly chose the flow times, contains the less-conducive-to-mixing value 2.  
Moreover, the best mixing was obtained by FSFT flows with $\tau=3,4$, which are time-periodic and hence most convenient in potential applications.

Randomizing the phase shift each time the direction of the flow switches seems to have a (slight but consistent) {\it negative effect} on the mixing rates for all $\tau\ge 3$, so mixing in these cases appears to be solely a result of stretching by the individual wedge flows.  Note that this effect is smaller for larger $\tau$, but in those cases there were also fewer additional random choices made during each run (besides the initial randomly chosen vertical and horizontal phase shifts).  While it need not be surprising that this randomness does not increase mixing, the strictly lower mixing rates in the RSXT cases vs.~their FSXT counterparts are unexpected.  We do not know what is the underlying reason, and further study of this phenomenon could be of interest.

We note that this could also suggest that the phase shift randomization proposed by Pierrehumbert in \cite{Pie,Pie2} might have a positive effect on mixing only when the flow direction switching is too frequent.  Moreover, in that case it might be more beneficial to increase flow times rather than randomize the phases, which could also be more easily implemented in real world situations.  Nevertheless, the sine flows in \cite{Pie,Pie2} have a different geometry from our wedge flows due to decreased stretching near the lines where their velocities are extremal (and hence their derivatives vanish), so simulations with these flows will be needed to determine whether the above conclusions also apply in this case.

The effect of randomizing the flow times is more difficult to discern.  Comparing the corresponding XSRT and XSFT cases suggests that it is negligible for $\tau\ge 7$, but this may not be unexpected since the variation in the flow times is small relative to the mean flow time $\tau$.  For $\tau=4,6$, one can observe some decrease of mixing rates when the flow times are randomized, but this reverses for $\tau=5$ 
(the XSRT cases with $\tau=3$ are again special because $2\in \{\tau-1,\tau,\tau+1\}$).  However, in the cases $\tau=4,5$ these differences can be explained at least in part by noticeably slower mixing when the (fixed) flow time is $\tau=5$ vs.~$\tau=4$, meaning that randomizing the flow time should slow down mixing for $\tau=4$ but speed it up for $\tau=5$.  In fact, in order to exclude  effects of both values 2 and 5, we also made 100 runs of each of the XSRT cases with flow times uniformly distributed in $\{3,4\}$; the obtained average mixing rates were much closer to the corresponding XSFT cases with $\tau=3,4$, namely $0.3905$ in the FSRT case and $0.3778$ in the RSRT case. Hence no clear pattern seems to emerge here.


\begin{figure}[htbp]
\centering
\begin{tabular}{|c|c|c|c|c|}
	\hline
	$\tau$ & FSFT & RSFT & FSRT & RSRT \\\hline\hline
	2 & --- & 0.0234 & 0.0569 & 0.0508 \\\hline
	3 & 0.0115 & 0.0168 & 0.0301 & 0.0250 \\\hline
	4 & 0.0136 & 0.0170 & 0.0144 & 0.0175 \\\hline
	5 & 0.0096 & 0.0161 & 0.0189 & 0.0193 \\\hline
	6 & 0.0105 & 0.0135 & 0.0146 & 0.0155 \\\hline
	7 & 0.0105 & 0.0113 & 0.0135 & 0.0168 \\\hline
	8 & 0.0102 & 0.0101 & 0.0126 & 0.0149 \\\hline
	9 & 0.0123 & 0.0139 & 0.0124 & 0.0154 \\\hline
	10 & 0.0125 & 0.0129 & 0.0139 & 0.0146 \\\hline
\end{tabular}
\caption{Standard deviations of (base-2) exponential mixing rates.}
\label{stdev}
\end{figure}

Let us also discuss variations in the data that yielded the average rates in Figure \ref{main_data}, since one may  wonder whether individual runs of our simulation exhibited mixing rates that are close enough to these averages.  
Figure \ref{stdev} shows the standard deviations of the 100 individual base-2 mixing rates in each case, demonstrating that the averages in Figure \ref{main_data} are also good approximations of the mixing rates of most individual runs.
Indeed, the standard deviations were no more than $0.0193$ in all cases with $\tau\ge 3$ except for the cases FSRT and RSRT with $\tau=3$, where the variation of mixing rates was  magnified due to  $2\in \{\tau-1,\tau,\tau+1\}$.  For instance, the largest difference between the mixing rate for an individual run with $\tau\ge 3$ and the  average mixing rate in the corresponding flow case  was $0.1020$ in the FSRT case with $\tau=3$; that run included a number of consecutive flow times 2 (and fixed phase shifts), a setup that leads to much slower mixing (see the next subsection).  


Finally, for completeness, we list the means and standard deviations of $T_{10}$ in Figure \ref{Tvalues}.


\begin{figure}[htbp]
\centering
\begin{subfigure}{0.5\textwidth}
\centering
\begin{tabular}{|c|c|c|c|c|}
	\hline
	$\tau$ & FSFT & RSFT & FSRT & RSRT \\\hline\hline
	2 & --- & 48.13 & 72.73 & 75.86 \\\hline
	3 & 27.58 & 28.88 & 30.09 & 30.80 \\\hline
	4 & 26.73 & 28.00 & 28.26 & 28.77 \\\hline
	5 & 28.22 & 29.48 & 28.56 & 29.46 \\\hline
	6 & 29.41 & 30.46 & 30.52 & 31.07 \\\hline
	7 & 31.28 & 32.00 & 31.79 & 31.86 \\\hline
	8 & 33.89 & 34.23 & 34.01 & 34.30 \\\hline
	9 & 34.90 & 35.93 & 34.92 & 35.92 \\\hline
	10 & 35.73 & 36.49 & 36.12 & 36.88 \\\hline
\end{tabular}
\end{subfigure}%
\begin{subfigure}{0.5\textwidth}
\centering
\begin{tabular}{|c|c|c|c|c|}
	\hline
	$\tau$ & FSFT & RSFT & FSRT & RSRT \\\hline\hline
	2 & --- & 3.67 & 18.95 & 17.66 \\\hline
	3 & 0.64 & 0.83 & 2.17 & 1.87 \\\hline
	4 & 0.49 & 1.11 & 0.84 & 1.05 \\\hline
	5 & 1.05 & 1.08 & 1.13 & 1.08 \\\hline
	6 & 1.54 & 1.18 & 1.37 & 1.30 \\\hline
	7 & 1.19 & 1.03 & 1.11 & 1.21 \\\hline
	8 & 0.65 & 0.45 & 1.04 & 1.21 \\\hline
	9 & 1.47 & 1.57 & 1.19 & 1.55 \\\hline
	10 & 1.32 & 1.66 & 1.49 & 1.68 \\\hline
\end{tabular}
\end{subfigure}

\caption{Means (left) and standard deviations (right) of $T_{10}$.}
\label{Tvalues}
\end{figure}

\subsection{Fixed Structures and Grid Scale Effects}
\label{S3.2}

Let us now turn to flow times $\tau=1,2$, starting with $\tau=2$.  The tables above are all missing data in the case FSFT with $\tau=2$, which is because one can easily show that these flows do not produce exponential mixing (hence we did not run our simulations in this case).  The segments $\{(s,s+\frac 34)\,|\,s\in (0,\frac 14)\}$ and $\{(s,s+\frac 14)\,|\,s\in (\frac 14,\frac 12)\}$ are fixed by $(\calV_{0}^{2} \circ \calH_{0}^{2})^2$, and mapped onto each other by $\calV_{0}^{2} \circ \calH_{0}^{2}$.  The same is true about the segments $\{(s,\frac 34-s)\,|\,s\in (\frac 12,\frac 34)\}$ and $\{(s,\frac 54-s)\,|\,s\in (\frac 34,1)\}$.  Obviously, $(\omega',\omega)$-shifts of these segments have the same relationship to $\calV_{\omega'}^{2} \circ \calH_{\omega}^{2}$.  Moreover, locally at any point on these segments, $(\calV_{\omega'}^{2} \circ \calH_{\omega}^{2})^2$ is represented by the matrix 
\[
\begin{bmatrix}
            -3 & 4 \\
            -4 & 5 \\
\end{bmatrix},
\]
which is similar to the $2\times 2$ Jordan block with diagonal elements $1$ (the eigenspace for this eigenvalue is generated by $(1,1)$).  It is easy to see from this that one can  at best hope for algebraic-in-time mixing in this case.  Figure \ref{mixing_example} in fact shows how mixing is inhibited near the above segments, as well as that it is faster elsewhere in the domain.

These structures of course do not survive the randomization in the RSFT, FSRT, and RSRT cases with $\tau=2$, but they show that flow time $\tau=2$ does not provide enough stretching and layering for the most efficient mixing.  This is compounded by flow times 1 occurring in the FSRT and RSRT cases, as data in Figure \ref{main_data} demonstrates (while the RSFT entry does show a decent mixing rate in that case, it is still well below the cases with $\tau\ge 3$).

The cases with $\tau=1$ were even slower mixers and we therefore did not perform their full simulations.  In the FSFT case, one can again easily identify a fixed structure of the mapping $\calV_{0}^{1} \circ \calH_{0}^{1}$.  Indeed the segments $\{(s,s+\frac 12)\,|\,s\in (0,\frac 12)\}$, $\{(\frac 12,s)\,|\,s\in (0,\frac 12)\}$, and $\{(s,\frac 12)\,|\,s\in (\frac 12,1)\}$ form a 3-cycle for this mapping and each is fixed by $(\calV_{0}^{1} \circ \calH_{0}^{1})^3$.  The latter mapping is locally represented by matrices
\[
\begin{bmatrix}
            3 & -2 \\
            2 & -1 \\
\end{bmatrix}
\qquad \text{and} \qquad 
\begin{bmatrix}
            -1 & 2 \\
            -2 & 3 \\
\end{bmatrix}
\]
on the two sides of the segment $\{(s,s+\frac 12)\,|\,s\in (0,\frac 12)\}$, and both are again similar to the $2\times 2$ Jordan block with diagonal elements $1$ (the eigenspace in both cases is again generated by $(1,1)$).  Thus there is no exponential mixing in this case either.

\begin{figure}[htbp]
    \centering
        \includegraphics[scale = 0.55]{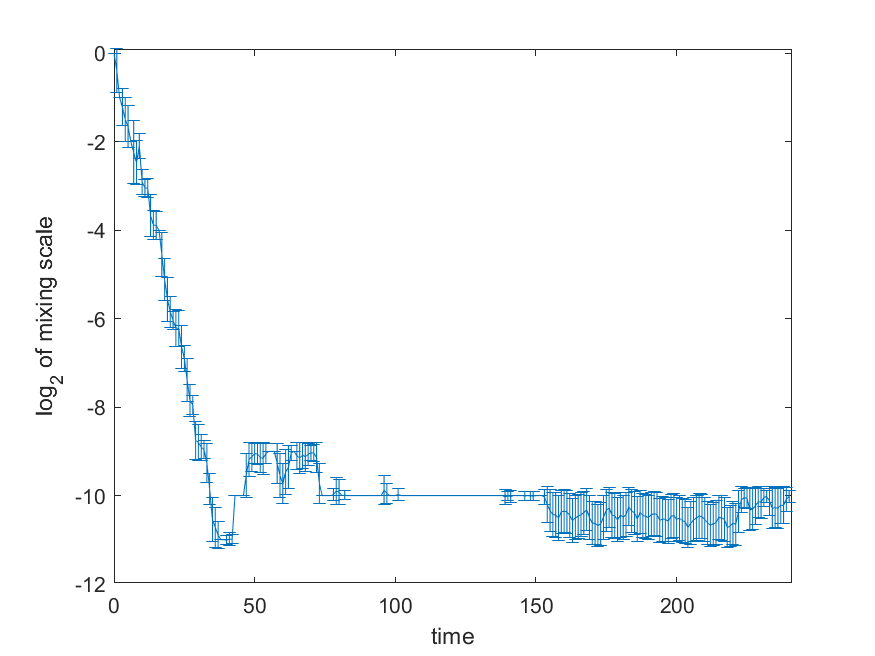}
        \includegraphics[scale = 0.55]{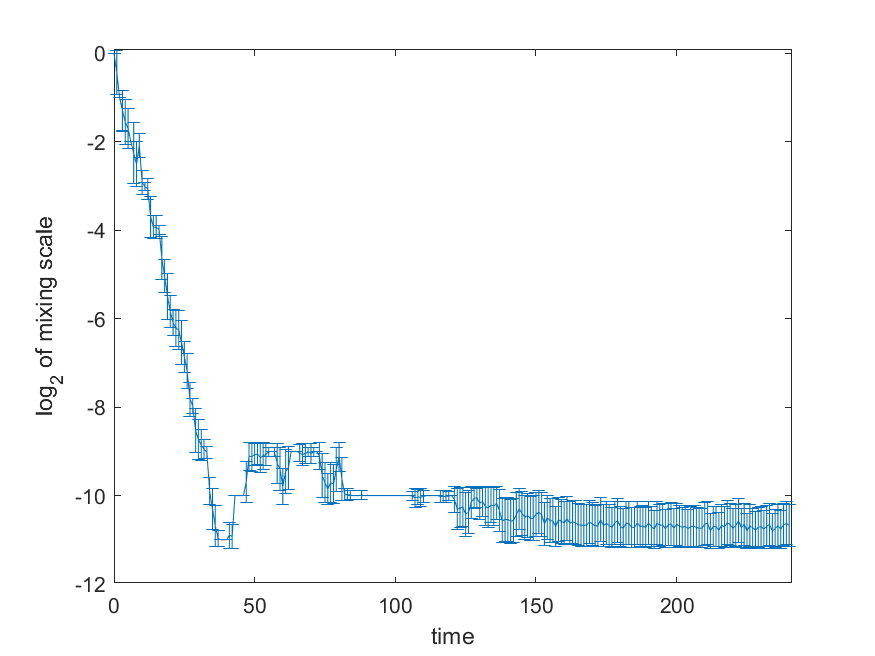}      
    \caption{Means and standard deviation error bars of the binary logarithms of the mixing scales for the 100 simulations 
   with $N = 15$ in the FSFT  (left) and RSFT   (right) cases with $\tau=8$, over an extended time interval.}
    \label{kinks}
\end{figure}

Finally, we mention here a curious phenomenon that we observed for the FSFT and RSFT cases with $\tau=8$.  In both cases, after the mixing scale $\sim 2^{-11}$ was reached and the observed exponential decrease stopped, the computed scale rebounded to $\sim 2^{-9}$ for some time (as if some {\it unmixing} were happening there).  This was then followed by a time interval (fairly long one, particularly in the FSFT case) where the mixing scale equaled $2^{-10}$ for {\it each} of the 100 simulations, before it stopped having this surprisingly uniform behavior and settled into a slightly more varied dynamic with values near $2^{-10}$ and $2^{-11}$.  
Figure \ref{kinks} contains means and standard deviation error bars of the binary logarithms $-n_k$ of the mixing scales   in these two flow cases.
We ran our simulations with a larger grid scale as well, and this behavior persisted in that case, albeit clearly at larger mixing scales.

These observations are of course completely irrelevant to the mixing theory, since they involve behavior on time intervals where the simulation cannot anymore capture the mixing dynamic due to the mixing scales being too close to the grid scale $2^{-15}$.  Nevertheless, we did not observe it or something similar for other values of $\tau$ (as well as in the XSRT cases where $\tau$ varies over time), so it may point to some special feature of the discrete grid dynamic for flow time $\tau=8$.  At the same time, the fact that this behavior is observed in the RSFT case as well, where the randomness in phase shifts would destroy any potential special structures present in the FSFT dynamic, is quite curious.  We do not currently have a candidate for the possible reason behind this phenomenon, and do not know why $\tau=8$ is the only flow time out of those we studied for which it occurs.

\Addresses


\begin{thebibliography}{99}
\bibitem{ACM} 
G. Alberti, G. Crippa, and A.L. Mazzucato, 
{\it Exponential self-similar mixing and loss of regularity for continuity equations},
C. R. Acad. Sci. Paris, Ser. I {\bf 352} (2014), 901--906.

\bibitem{ACM2} 
G. Alberti, G. Crippa, and A.L. Mazzucato,
{\it Exponential self-similar mixing by incompressible flows}, 
J. Amer. Math. Soc.  {\bf 32} (2019), 445--490.





\bibitem{BBP} J. Bedrossian, A. Blumenthal, and S. Punshon-Smith, 
{\it Almost-sure exponential mixing of passive scalars by the stochastic Navier-Stokes equations},
Ann. Probab., to appear.

\bibitem{BOS} 
Y. Brenier, F. Otto, and C. Seis, 
{\it Upper bounds on the coarsening rates in demixing binary viscous fluids},
SIAM J. Math. Anal. {\bf 43} (2011), 114--134.

\bibitem{Bressan} 
A. Bressan,
{\it A lemma and a conjecture on the cost of rearrangements},
Rend. Sem. Mat. Univ. Padova {\bf 110} (2003), 97--102. 

\bibitem{B2} A. Bressan. Prize offered for the solution of a problem on mixing flows. http://www.math.psu.edu/bressan/PSPDF/prize1.pdf, 2006.




\bibitem{CL} 
G. Crippa and C. De Lellis, 
{\it Estimates and regularity results for the DiPerna-Lions flow},
J. Reine Angew. Math. {\bf 616} (2008), 15--46.

\bibitem{CLS} 
G. Crippa, R. Luc\'{a}, and C. Schulze,
{\it  Polynomial mixing under a certain stationary Euler flow},
Phys. D {\bf 394} (2019), 44--55.



%
%
%

\bibitem{ElgZla} 
T.M. Elgindi and A. Zlato\v s, 
{\it Universal mixers in all dimensions}, 
Adv. Math. {\bf 356} (2019), 106807, 33 pp. 

\bibitem{IKX} 
G. Iyer, A. Kiselev, and X. Xu, 
{\it Lower bounds on the mix norm of passive scalars advected by incompressible enstrophy-constrained flows},
Nonlinearity {\bf 27} (2014), 973--985.


\bibitem{LTD} 
Z. Lin, J. L. Thiffeault, and C. R. Doering,
{\it Optimal stirring strategies for passive scalar mixing},
J. Fluid Mech. {\bf 675}, 465--476.


\bibitem{LLNMD} 
E. Lunasin, Z. Lin, A. Novikov, A. Mazzucato, and C. R. Doering, 
{\it Optimal mixing and optimal stirring for fixed energy, fixed power, or fixed palenstrophy flows},
J. Math. Phys. {\bf 53} (2012), 115611, 15pp.


%

\bibitem{OSS} 
F. Otto, C. Seis, and D. Slep\v{c}ev,
{\it Crossover of the coarsening rates in demixing of binary viscous liquids},
Commun. Math. Sci. {\bf 11} (2013), 441--464.

\bibitem{Pie} 
R. Pierrehumbert, 
\it Tracer microstructure in the large-eddy dominated regime, 
\rm Chaos, Solitons \& Fractals {\bf 4} (1994), 1091--1110.

\bibitem{Pie2} 
R. Pierrehumbert, 
\it Lattice models of advection-diffusion, 
\rm Chaos {\bf 10} (2000), 61--74.

\bibitem{S} 
C. Seis,
{\it  Maximal mixing by incompressible fluid flows},
Nonlinearity {\bf 26} (2013), 3279--3289.

\bibitem{S2} 
D. Slep\v{c}ev,
{\it Coarsening in nonlocal interfacial systems}, 
SIAM J. Math. Anal. {\bf 40} (2008), 1029--1048.




\bibitem{T} 
J.-L. Thiffeault,
{\it Using multiscale norms to quantify mixing and transport},
Nonlinearity {\bf 25} (2012), 1--44.


\bibitem{YaoZla} 
Y. Yao and A. Zlato\v{s},
{\it  Mixing and un-mixing by incompressible flows},
J. Eur. Math. Soc. {\bf 19} (2017), 1911--1948. 


\end{thebibliography}
\end{document}